\documentclass{article}
\usepackage{maa-monthly}

\theoremstyle{theorem}

\newtheorem{problem}{Problem}
\newtheorem{particularproblem}{Particular case of Problem}
 
\theoremstyle{definition}
\newtheorem*{definition}{Definition}

\begin{document}

\title{A list of problems in Plane Geometry with simple statement that remain unsolved}
\author{L. Felipe Prieto-Mart\'inez}

\maketitle

\begin{abstract}
 This article contains a short and entertaining list of unsolved problems in Plane Geometry. Their statement may seem naive and can be understood at an elementary level.  But their  solutions have refused to appear for forty years in the best case.
\end{abstract}

\section*{Introduction}

It is possible to find several works in the bibliography and many web pages on the internet dedicated to open problems with an elementary statement. This is an attractive topic for various reasons and for a broad audience. To sum up, we could highlight two of them: (i) if anybody with an elementary education can understand the statement, anybody with an elementary education can think about its mysteries and (ii) if anybody with an elementary education can understand the statement of a problem and it remains unsolved, then it is because it requires revolutionary techniques or viewpoints which have not been developed yet. Certainly, their statements may be elementary, but their solutions are not.

Not every branch of Mathematics is equally accessible to a non-specialist readers. Two of the ones that are closest to all kind of  public are Combinatorics and Geometry (and so their intersection, Geometric Combinatorics). Geometry has the additional advantage of involving much more drawing. Good examples of compilations of problems of the kind we are interested in are, for instance, \cite{BMP, CFG, G, KW, W}. More specifically, problems in  Plane Geometry seem to be specially captivating and less demanding with respect to the graphic skills of the reader.

In this article, a list of  \emph{Problems in Plane Geometry with simple statement that remain unsolved} is offered.  Due to the abundance of similar works, and to offer the reader a quality product in these pages,  this list has been made according to some not very ambitious criteria. For a problem to be included in the list, it must  satisfy all the following necessary conditions:

\begin{enumerate} 

\item[(1)] It must be ``consolidated'' (more than 40 years from its publication and some attention from the mathematical community).

\item[(2)] It should be understood by an average final year undergraduate student in Mathematics, and its statement by even younger persons (additional information about the concepts involved has been included when necessary). 

\item[(3)] Including this additional information, its statement should be ``short''.

\item[(4)] Its statement must deal with very concrete objects, rather than with general objects.

\end{enumerate}

\noindent Surely, the list of six problems presented does not cover the family of problems matching this criteria, so the reader may feel that it is rather arbitrary. Certainly, it is not possible to coin such a list. The problems have been sorted from the newest to the oldest.  In the last section some other problems that ``almost fit'' the conditions are briefly commented. When a problem resists to be solved, variations of the statement are studied. Information about these variations  has also been intentionally excluded.

The criteria followed make the list to include problems that are  quite different in flavour: some of them belong to the realm of Geometric Combinatorics,  some others to Optimal Dessign of curves or bodies (see next section) when restricted to some geometric condition and the last one is a problem of existence.

The motivation to write this article is to provide a list of inspiring questions interesting for a broad audience (including secondary education students) in a concise style. Each reader may read this work with a different level of depth. The author has used this notes successfully to introduce real research and some concepts of higher Mathematics to his students.  It should not be regarded as a survey of the problems mentioned made by an expert. It just contain the basic concepts and bibliography (compiling some of the attempted solutions). Thus, we hope that it turns out to be sufficiently easy to read  that it stimulates the reader curiosity, giving his or her the option of further deepening into their own question.

\section*{Formalization of some intuitive geometric concepts}

As any mathematician knows, if a concept in Geometry is intuitive, it is not necessarily easy to formalize. Let us recall that there is a quote attributed to Felix Klein (see \cite{BOYER}) that states that \emph{everyone knows what a curve is, until he has studied enough mathematics to become confused through the countless number of possible exceptions}.

Despite the naive nature of this work, it is important not to lose rigor, typical of Mathematics. This is a difficulty that other works have overcome (see for example \cite{BMP}). And that needs to be done carefully. We want this article to contain, at least the basic concepts, although we cannot include topological and geometrical detail.

\begin{definition} A \textbf{parametrized plane curve} is a set $C\subset\mathbb R^2$ together with a surjective continuous funcion $\gamma:[a,b]\to\mathbb C^2$. $\gamma(a)$ is the so called initial point and $\gamma(b)$ is the final point.

\end{definition}

\begin{definition}[following \cite{GERBER}] A \textbf{plane region} is any closed connected subset $R\subset \mathbb R^2$.

\end{definition}

\begin{definition}[following \cite{BMP}] A \textbf{plane body} is a set in the plane homeomorphic to the closed unit disk. We say that two bodies do not overlap if they do not have an interior point in common. And we say that two of them touch if they do not overlap but they have a common boundary point.

\end{definition}

\section{Kobon's Triangle Problem}

This problem is considered to fist appear, with its definitive statement, in the 1978 book of  puzzles \cite{F}. The author is Kobon Fujimura. It belongs to the branch of mathematics known as \emph{Combinatorial Geometry}.

\begin{figure}
\includegraphics{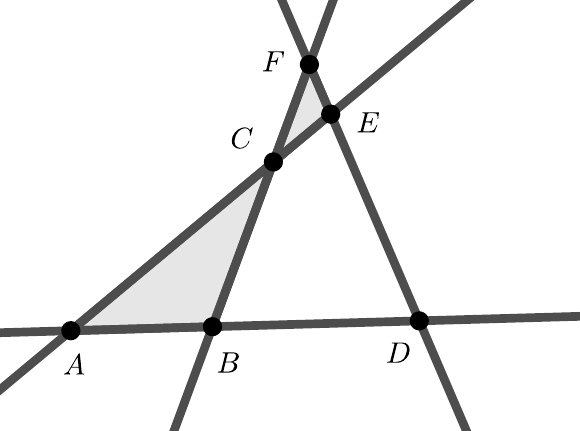}
\caption{Unlike other problems, we are not counting the triangle $ADE$ because it is cut by the line through $B$, $C$, $F$. With this agreement, we can avoid the problem of overlapping triangles ($ADE$ would overlap the triangle $ABC$).}
\label{kobon.explanation}
\end{figure}

\begin{problem} Determine the largest number $K(n)$ of (nonoverlapping, uncut, see Figure \ref{kobon.explanation}) triangles that  can be drawn using only $n$ straight lines in the plane.

\end{problem}

Obviously, it is not possible to draw any triangle with 1 or 2 lines. We can obtain 1 triangle with 3 lines. The reader may easily find a configuration with 4 lines and 2 triangles and another one with 5 lines and 5 triangles.

The sequence $(K(1),K(2),K(3),\ldots)$, called \emph{Kobon's Sequence}, appears in the \emph{Online Encyclopedia of Integer Sequences} \cite{OEIS} with the code $A006066$. Some of its terms are know, and can be found there. The terms $K(1),\ldots, K(9)$ are:
$$0,0,1,2,5,7,11,15,21 $$

\begin{figure}
\includegraphics[width=10cm]{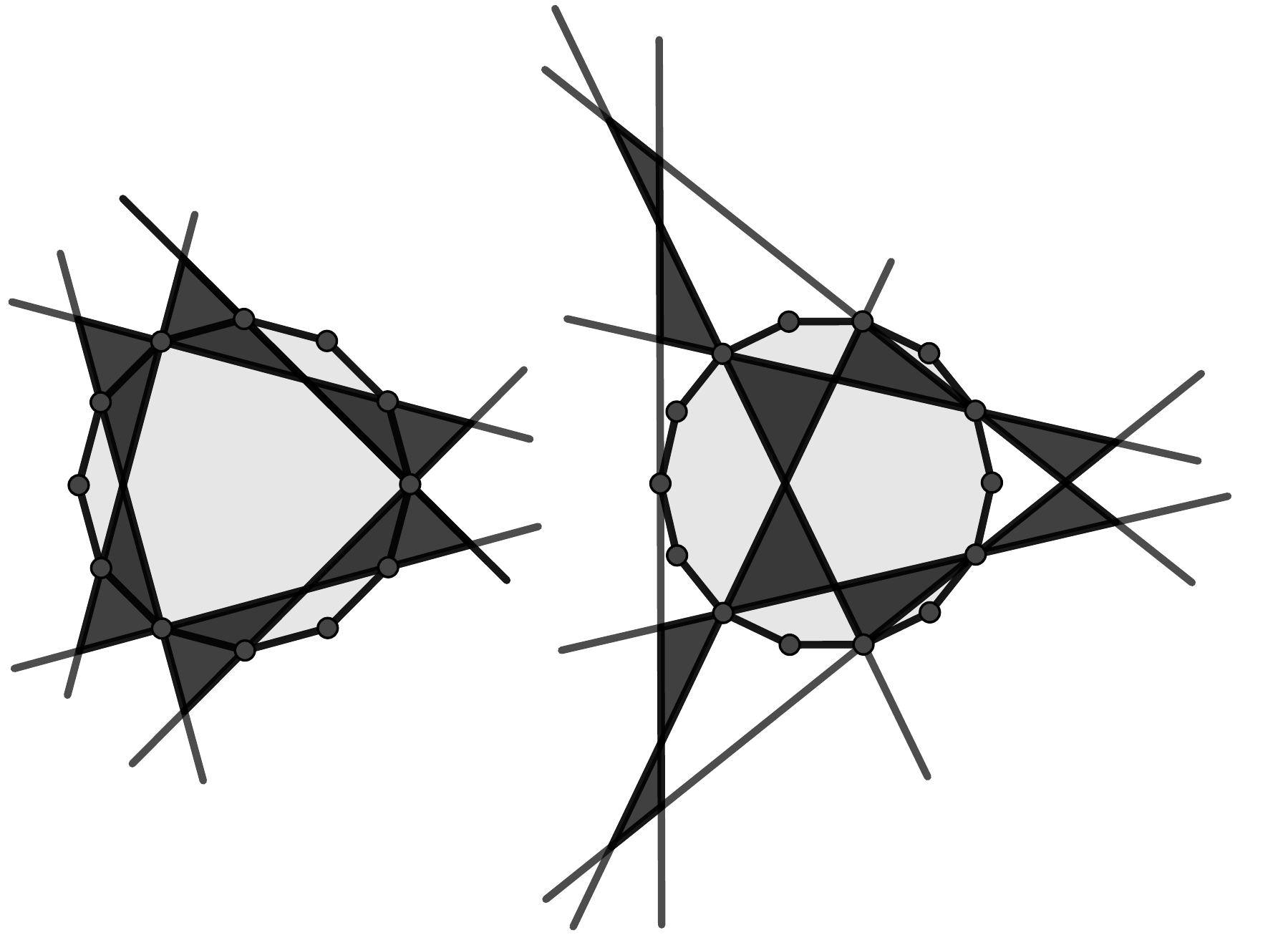}
\caption{The configuration by F\"uredi and Pal\'asti for $n$ lines consists in an adequate  choice of the lines containing the diagonals of a regular polygon of $2n$ sides. In this picture the corresponding configuration for $n=6,7$ is showed. None of them is optimal, since $K(6)=7, K(7)=11$.}
\label{kobon.furedipalasti}
\end{figure}

Each time we find a configuration with $n$ lines and $m$ triangles we are proving that $K(n)\geq m$.  There exist some general constructions that offer an arrangement of lines with a large number of triangles for each $n$. So, from this constructions, we can obtain lower bounds for the numbers $K(n)$. One of the best general constructions known is the one by F\"uredi and Pal\'asti in \cite{FP} (see Figure \ref{kobon.furedipalasti}):
$$\frac{n(n-3)}{3}\leq K(n) $$

But a difficult obstacle is that it is not easy to prove that a given configuration is optimal. One possible strategy, which has been used successfully for some values of $n$, is to find some upper bound $M(n)$ for the number $K(n)$ and then to find a configuration with $n$ lines and $M(n)$ triangles. One of the results that have made this way of action possible is the one by Cl\'ement and Bader in \cite{CB}, who found an acceptable upper bound for the numbers $K(n)$:
$$K(n)\leq \begin{cases} \lfloor \frac{n(n-2)}{3}\rfloor& \text{if }n\equiv 0,2\mod 6\\ \lfloor \frac{n(n-2)}{3}\rfloor & \text{in other case}\end{cases}  $$

\noindent Unfortunately, this strategy is not always satisfactory: frequently there is a gap between the best known configuration and the best known upper bound. This is the case, for example, of $n=10$:

\begin{particularproblem} {The best configuration known with 10 lines has 25 triangles. But the bound by Clement and Bader only guarantees that $K(10)\leq 26$. Maybe the bound is not tight enough? or it is possible to find a configuration with 10 lines and 26 triangles?}
\end{particularproblem}

Let us recall that some values of $K(n)$ are known for $n>10$. The reader may find an updated survey concerning the state-of-the-art of this problem in \cite{MP}. This article also includes most of the relevant configuration for small values of $n$ (like the best one known for $n=10$).

\section{Moving Sofa Problem}

This problem is attributed (in its final statement) to L. Moser, and appears in his paper \cite{M} published in 1966. Informally, the problem can be stated as follows

\begin{problem} Determine the plane region $S$ (from now on and in this section we will call it \emph{sofa}) of largest area which can be ``moved around a right-angled corner in a hallway of width one'' (see below).

\end{problem}

We impose that the sofa cannot be rotated in the space, it can only be dragged across the floor. So we are considering this as a problem in the plane. So, formally, as explained in \cite{GERBER}, we take the ``hallway'' to be the plane region:
$$H=\{(x,y):x,y\leq 1, \text{ either $x$ or $y$ are greater or equal to 0} \}$$ 

\begin{figure}
\includegraphics[width=4cm]{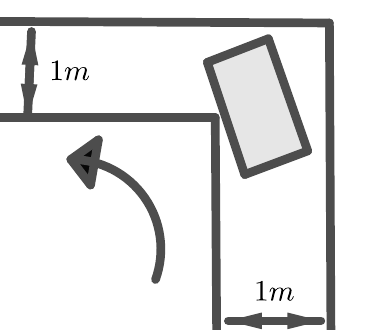}
\caption{Illustration of the Moving Sofa Problem.}
\label{sofa.intro}
\end{figure}

\noindent $S$ must be a subset of $H$ such that if $(x,y)\in S$ then $0\leq x\leq 1$, $y\leq 1$. We must be able to  ``rigidly move $S$ through $H$'. This is a very reasonable sentence (one may understand a situation as the one in Figure \ref{sofa.intro}) which is not easy to write with rigor. Formally, for a parameter $t$ representing time and varying in a closed interval $I$, we need  to find a uniparametric family of rigid motions in the plane $\Phi_t$ such that (1) for every $t\in I$ ``the sofa is inside the hallway'' ($\Phi_t(S)\subset H$) and (2)  ``the sofa moves continously and without teleportation through the corridor'' (for any fixed $(x_0,y_0)\in S$, $\gamma(t)=\Phi_t(x_0,y_0)$ is a parametrized plane curve).

\begin{figure}
\includegraphics[width=6cm]{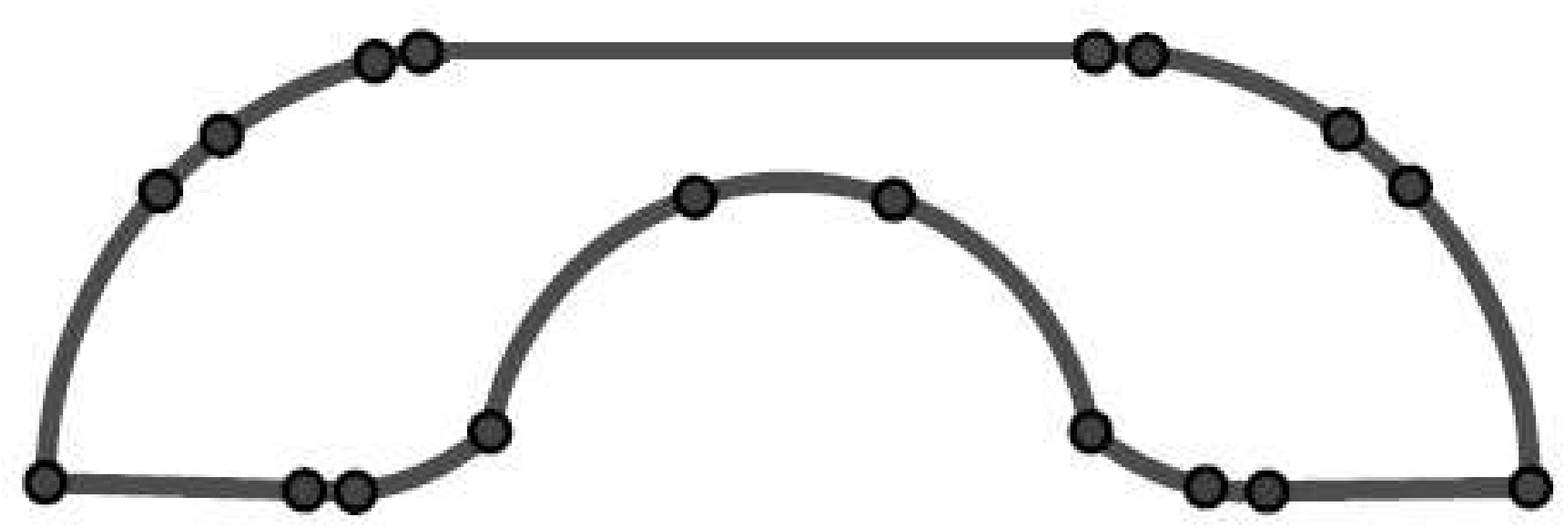}
\caption{Gerver's sofa is obtained joining 18 arcs.}
\label{sofa.gerber}
\end{figure}

\begin{figure}
\includegraphics[width=6cm]{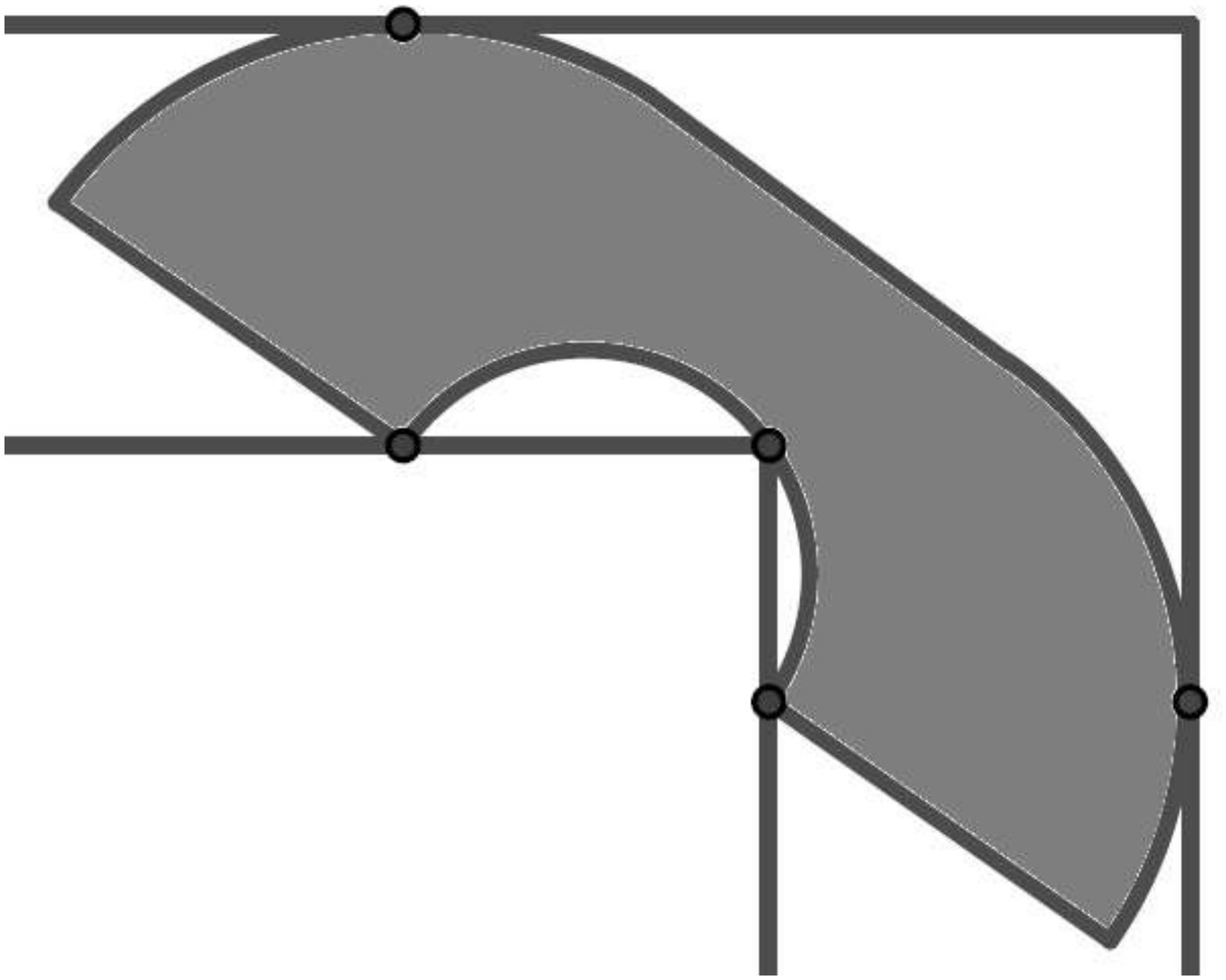}
\caption{This is Hammersley's sofa, very similar to Gerver's sofa, moving through the corner.}
\label{sofa.hammer}
\end{figure}

The sofa needs not to be convex. For instance, the shape of the better known solution (its area is approximately 2.2195) reminds to an old telephone (see Figure \ref{sofa.gerber}). It was proposed by  J. Gerver (see  \cite{GERBER} for a precise description) and it is conjectured to be the optimal sofa. To see how such an ``old telephone'' is moved around the corner see Figure \ref{sofa.hammer}.

In the bibliography we may find two different types of attempts  to obtain the optimal sofa. On the one hand, a computational approach  has  been considered. See, for instance, the article by P. Gibbs \cite{GIBBS}. The solution obtained in this work is indistinguishable from the one by Gerver. On the other hand, also analytical methods have been used. See the paper by D. Romik \cite{R.DIFF}.

Finally we would like to highlight the work by Y. Kallus and  D. Romik \cite{R.ADV}. In this article, it is proved that the area of the optimal sofa $A$ is less or equal to 2.37 (Gerver's sofa proves that $A\geq 2.2195$).

We recommend the enthusiastic introduction of \cite{R.DIFF} to see the impact of this problem in the mathematical community and a survey with a lot of examples and the video \cite{R} of this same author explaining the problem and some of its features in a simple way.

\section{Moser's worm problem}

This problem was published in \cite{M2} also by L. Moser and also in 1966.


\begin{problem} Define a \emph{worm} to be any plane curve of length one. Define a \emph{cage} to be a plane region that can \emph{accommodate} any worm, that is, the worms have a fixed shape that cannot be altered, but they can be rotated and translated  to fix into the region. Find, if it exists, the cage of smallest area.

\end{problem}

Lately, most of the authors listed below have directed their attention to the case in which the cage is additionally required to be convex. We refer to this problem as  Moser's convex worm problem.

It is easy to find simple regions that can accommodate any worm. For instance, a disk of diameter 1 is a cage. Gerriets and Poole in \cite{GP} proved that the rhombus with diagonals of lengths 1 and $1/\sqrt{3}$ (and thus obtained ``gluing together'' two equilateral triangles of sidelength equal to $1/\sqrt{3}$) is also a cage. We recommend the reader to check this proof. Quite recently, C. Panraksa and W. Wichiramala showed that the $30º$ circular sector is also a cage. This is the best convex cage known, for this moment. But for the general case, the best cage known, up to the author knowledge, is the one explaines in \cite{NP} which is non-convex, with an area of 0,2604 approximately.

In a similar fashion to what happened for Problem 1, one of the main obstacles is that there is no known general strategy for showing that a given planar region is a cage. Gerriets and Poole tried to find such a strategy: they conjectured that a region is a cage if and only if it  can accommodate any \emph{polygonal worm}, that is, any \emph{polygonal chain}  (a connected and finite sequence of segments, where the starting point of each element is the endpoint of the previous one) with three elements and of length 1. But this was later disproved by  Panraksa, Wetzel and Wichiramala in \cite{PWW}. They showed, moreover, that ``three'' cannot be replaced by any other natural number in the previous statement to obtain a valid result.

As suggested in the statement, the following situation could also be possible: it may exists a value $A$, such that we can find cages with area as close as we want to $A$ but such that there is no cage of area $A$. It can be  that his cannot happen for the Moser's convex worm problem. It is a consequence of the Blaschke selection theorem. See the reference to this fact made in \cite{NPL}.

\section{Erd\"os-Oler packaging problem}


We will first state the general problem:

\begin{problem} \label{problem.pack} For each $n\in\mathbb N$, find the smallest triangle in which you can \emph{package} $n$ circles of radius 1. By \emph{package} we mean that the $n$ circles must be completely contained in the the triangle without overlapping.

\end{problem}

This problem is an abstraction of a real life question. It belongs to the larger family of circle packaging problems. The word ``triangle'' in the statement may be replaced by``circle'', ''square'' and``polygon congruent to one given'' to obtain another problem as interesting as this one. It is not easy to decide which is the first mention of this particular statement. But we will provide concrete information and authorship of the particular case explained below.

\begin{figure}
  \centering
    \includegraphics[width=50mm]{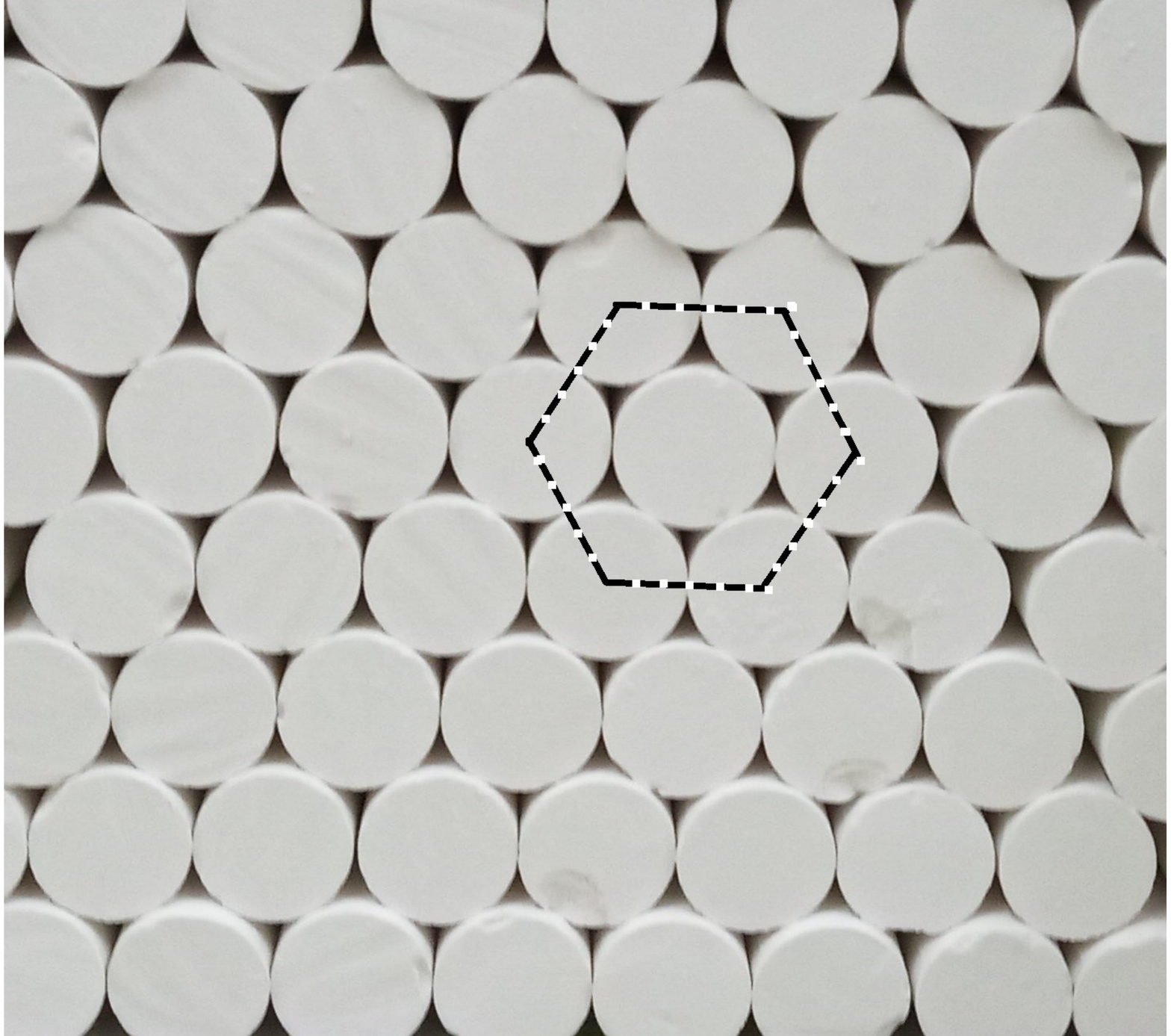}
\caption{Thue's Theorem states that the densest way to store circles in the plane is this hexagonal pattern exhibited in this box of chalks. Despite the name of the Theorem, it seems that it was J. P. Lagrange the first one to prove it.}
  \label{im.tizas}
\end{figure}

It is also related to the problem of storing circles in the plane in the densest way and to Thue's Theorem (see Figure \ref{im.tizas}).

\begin{figure}
  \centering
    \includegraphics[width=50mm]{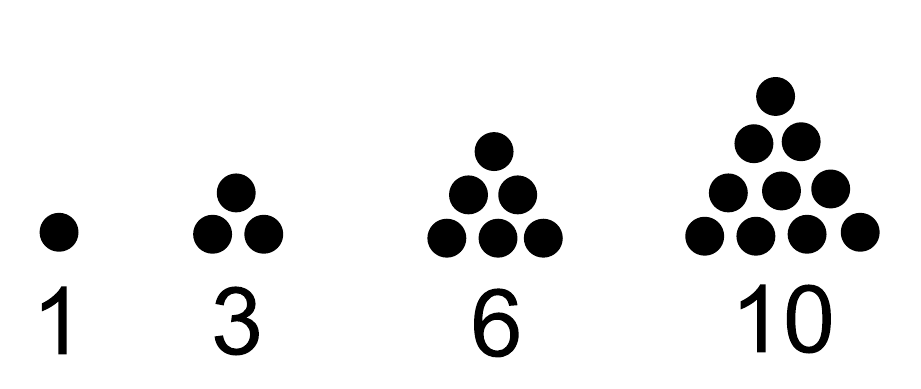}
\caption{The first triangular numbers and their geometric pattern.}
  \label{im.triangular}
\end{figure}

\begin{figure}
  \centering
    \includegraphics[width=40mm]{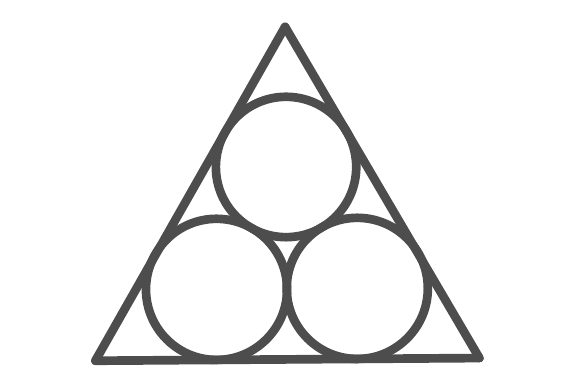}
 \caption{Here we may find the optimal packaging of three circles.}
 \label{pack3}
\end{figure}

Recall at this point that the sequence of \textbf{triangular numbers} (code 000217 in \cite{OEIS}) is $1, 3,6,10,\ldots$ Triangular numbers correspond to the geometric pattern exhibited in Figure  \ref{im.triangular} and so the sequence satisfies,  for $n\geq 2$, that the $n$-th element $a_n$ is obtained as $a_{n}=a_{n-1}+n$. Equivalently $\forall n\geq 1$, $a_n=(n(n+1)/2)$. This sequence is directly related to our question. Problem 4 was quickly solved for any triangular number $n$ (see \cite{O}). The way to organize the circles is the obvious one (see Figure \ref{pack3}). To verify this just see that, if the optimal way to store the circles is not this one, this would contradict Thue's Theorem.

 The triangles needed up to $n=18$ are well known (see \cite{MS}).  Nurmella in \cite{N} provided a numerical approach to this problem for the cases $n\leq 36$.

\begin{figure}
  \centering
    \includegraphics[width=80mm]{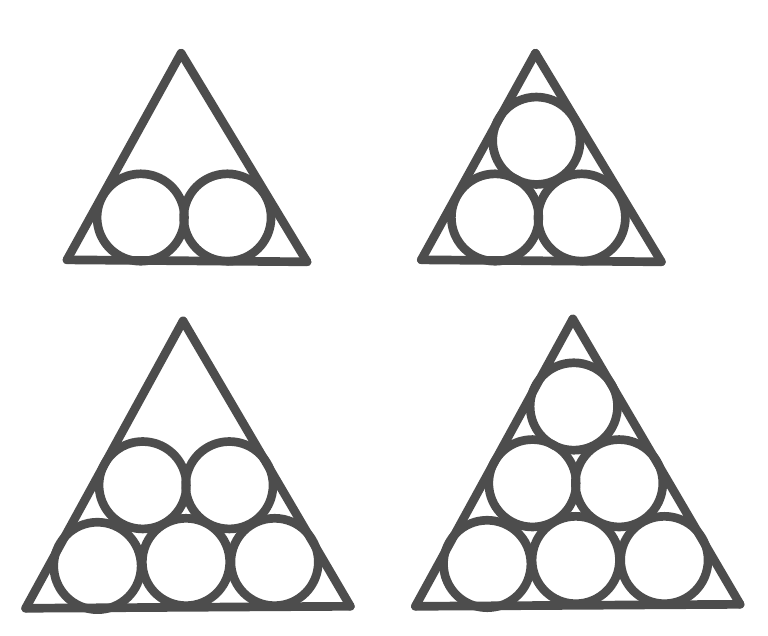}
\caption{Illustration of the first cases of packaging of circles for triangle numbers and for triangle numbers minus one.
}
 \label{pack.conjetura}
\end{figure}

Finally, we would like to remark that the following particular case of our problem is of great interest too. It was discussed by Erd\"os and Oler, and published in \cite{O} in 1961. See Figure \ref{pack.conjetura}.

\setcounter{particularproblem}{3}
\begin{particularproblem} {Let $n$ be a triangular number. The smallest equilateral triangle needed to package $n-1$ circles of radius 1 and to package $n$ circles of radius 1 is the same.}
\end{particularproblem}

\section{Bellman's lost in the forest problem}


The author of this problem is R. E. Bellman, who published it in his article \cite{B} in 1956. Suppose that a hiker is lost in a forest, that is, he knows the shape of the forest (imagine he has a map). (1) He does not know in which point of the forest he is placed, neither in which direction he is looking. (2) He wants to decide which is ``the best'' way to move in order to scape from the forest. From the point of view of the hiker, it is not clear what does it mean  ``the best''. This will be discussed below. A formal statement could be the following:

\begin{problem} Let $F$ be a plane region that we will call  \emph{forest}. Find a parametrized curve $C$ such that (1) for every $P\in F$, for every curve $\widetilde C$ congruent to $C$ and with initial point $P$, $\widetilde C$ ``scapes from $F$'' (it is not contained in $F$) and (2) among them find the ``best'' possible one.

\end{problem}

Here by ``best'' we will interpret \emph{the shortest time of escape in the worst case} or, in other words, \emph{the shortest trajectory that guarantees to escape}. But, in the bibliography, other authors interpret ``best'' as  \emph{the trajectory that minimizes the average time to escape}. In our case, the problem is closely related to Moser's worm problem: if any curve of length $1$ does not guarantee to escape from a forest, then this forest is a \emph{cage}.

The problem is still open for very simple choices of $F$. The curve requested may not be unique. A good survey about this topic is \cite{FW}. A summary, extracted from this article, of the situation is the following:

\begin{itemize}

\item Using what we know about the rhombus by Gerriets and Poole (described before), we can define a \emph{fat forest of diameter $L$} to be any forest $F$ such that (1) it contains two points $P,Q$ such that the distance between $P,Q$ coincides with the diameter of $F$ and (2) $F$ contains a rhombus with diagonals of lengths $L$ and $L/\sqrt{3}$ where the main diagonal is the segment $PQ$.Then we can easily prove that any fat forest of diameter $L$ needs a trajectory of length at least $L$, and that the best trajectory is, in fact,  a straight line segment of length exactly $L$. This includes $F$ being a disk, some cases of circular sectors and regular polygons for $n>3$.

\item For an infinite band, the solution is Zalgaller's curve which, surprisingly, the answer is not a straight line segment (see Figure \ref{B.Z}). This solution also applies for \emph{eccentric} rectangles.

\begin{figure}
  \centering
    \includegraphics[width=110mm]{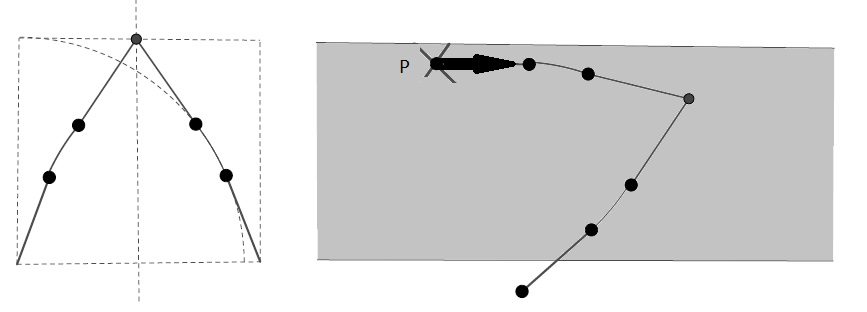}
 \caption{On the left, we can see Zalgaller's curve. On the right, we can see a diagram representing a hiker starting from the point $P$ initially looking in the direction of the arrow, and how Zalgaller's curve takes him out of the forest. It can be proved that no  shorter curve guarantees to escape from the forest.}

 \label{B.Z}
\end{figure}

\item The problems remains unsolved for triangles. There is a conjecture for the best curve for equilateral triangles in \cite{BES}.

\end{itemize}

\section{The inscribed square problem}

Our last question is my favourite one. It is due to O. Toeplitz, who proposed it in a conference in 1911 about problems in Topology. Recall that a \emph{Jordan curve} is a parametrized curve which is \emph{closed} (the initial and final points coincide) and \emph{simple} (injective, so ``it does not meet the same point more than once''). According to the Jordan curve theorem, Jordan curves divide the plane in two connected regions: one of them bounded, called the \emph{interior region}, and the other one unbounded and called the \emph{exterior region}.

\begin{figure}
  \centering
    \includegraphics[width=110mm]{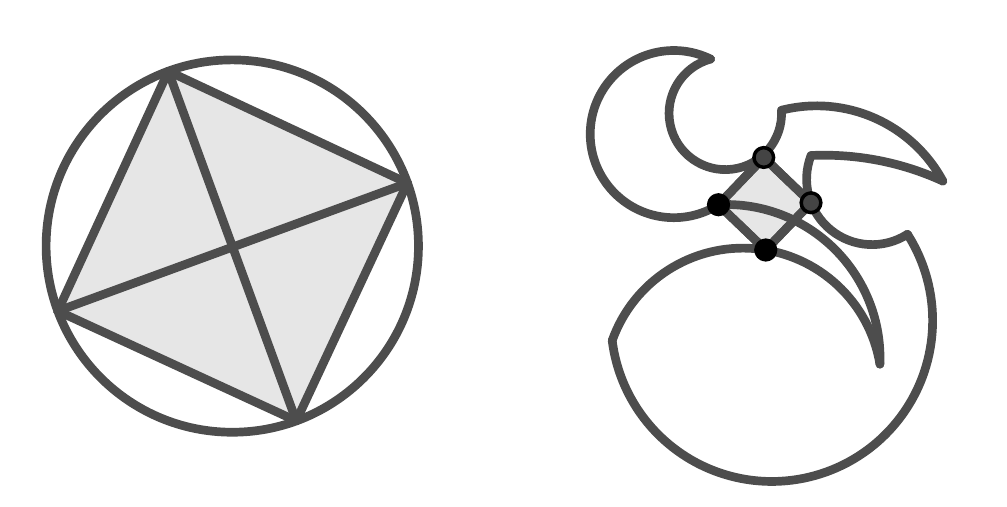}
 \caption{Circles are Jordan curves, and they inscribe infinitely many squares. On the left, we can see a Jordan curve slitghtly more complicated and a square inscribed in it.}

 \label{cuadradoinscrito}
\end{figure}

\begin{problem} Decide if any Jordan curve inscribes a square. Here ``inscribe'' means that the four vertices of the square are contained in the curve, not that the entire square is contained in the interior region (see Figure \ref{cuadradoinscrito}).

\end{problem}

Great mathematicians as T. Tao have got involved in this problem and we can find several attempts to solve it, very different in nature.

The main obstacle is that Jordan curves can be very complicated. For instance, they may not be \emph{rectifiable} (the may have infinite length).  For families of ``simple'' Jordan curves (convex curves, polygons, piecewise analytic curves, symmetric curves with respect to a line or point, locally monotone curves, $C^1$ curves,\ldots) the problem has been already solved. But the general question have remained open for more than a century.

Some results have also be obtained showing that any Jordan curve inscribes polygons of some given families (triangles congruent to one given, rhombi, rectangles,\ldots).

The survey \cite{MAT1} is strongly recommended, and it sums up  the current situation of the problem. See also \cite{MAT2} and the introduction of \cite{T}.

The main obstacle when trying to solve this problem is that the condition ``being a Jordan curve'' is purely \emph{topological}. We know almost nothing about the curve and so we have hardly any property to use. Simple geometrical strategies that we may dessign will always fail since there are a lot of monster curves in the family of Jordan curves.

\section*{Final comments}

Much more problems could have been included in the list. I have decided to, at least, name and briefly  comment some more of them and their cause of exclusion with respect to the conditions explained in the introduction.

\begin{itemize}

\item \textbf{Einstein problem (2, 3).} It concerns aperiodic tilings in the plane. See \cite{ST} for more information.

\item \textbf{Hammer X-ray problem (3).}  It deals with the number of ``ideal tomographic pictures'' (with no noise) needed for the identification or reconstruction of plane regions from information about their chords. See \cite{GARDNER, GM, GK}.

\item \textbf{Rado's covering problem (2,3).} A question about  families of squares in the plane with their sides parallel to the coordinate axis. See \cite{BDJ}.

\item \textbf{Other packaging problems (4).} A lot of variations of Problem \ref{problem.pack} are of great interest in Applied Mathematics. For instance packaging in any other regular polygon, or in circles, or involving circles of different sizes.

\item \textbf{$F$-chordal point problems (4).} The equichordal point problem was first proposed in 1916-17 and solved in 1997 by Rychlik in \cite{Ri}.  This problems asked about  existence and uniqueness of curves having two \emph{equichordal points}, that is, points such that every chord passing through it has the same length. A generalization of equichordal points are $F$-chordal points, which, for a given curve, satisfy  certain relations for the chords  passing through them. Several variations of the initial problem remain open. See \cite{P} for the basic definitions.

\item \textbf{Conjecture about the circle (1).} This problem could have been included in the previous family, but it deserves to be highlighted. Show that the circle is the unique $C^\infty$ Jordan curve with an equipower point (see \cite{Z} for the definition). It is easy to prove that \emph{lunes} actually have an equipower point, but they are only piecewise $C^\infty$.

\item \textbf{Heesch's problem (2,3).} For some region $D$ in the plane, determining the maximum number of layers of (congruent) copies of $D$ that can surround it. See \cite{MANN}.

\end{itemize}

I hope that you have enjoyed playing with this problems.



\vfill\eject

\end{document}